# CONSTRAINED BROWNIAN MOTION: FLUCTUATIONS AWAY FROM CIRCULAR AND PARABOLIC BARRIERS

### By Patrik L. Ferrari[1] and Herbert Spohn

### *Technische Universität München*

Motivated by the polynuclear growth model, we consider a Brownian bridge $b(t)$ with $b(\pm T) = 0$ conditioned to stay above the semicircle $c_T(t) = \sqrt{T^2 - t^2}$. In the limit of large $T$, the fluctuation scale of $b(t) - c_T(t)$ is $T^{1/3}$ and its time-correlation scale is $T^{2/3}$. We prove that, in the sense of weak convergence of path measures, the conditioned Brownian bridge, when properly rescaled, converges to a stationary diffusion process with a drift explicitly given in terms of Airy functions. The dependence on the reference point $t = \tau T$, $\tau \in (-1, 1)$, is only through the second derivative of $c_T(t)$ at $t = \tau T$. We also prove a corresponding result where instead of the semicircle the barrier is a parabola of height $T^\gamma$, $\gamma > 1/2$. The fluctuation scale is then $T^{(2-\gamma)/3}$. More general conditioning shapes are briefly discussed.

**1. Introduction and main results.** We consider the Brownian bridge $b(t)$ over the time interval $[-T, T]$, $T > 0$, $b(-T) = b(T) = 0$, conditioned to lie above the semicircle $c_T(t) = \sqrt{T^2 - t^2}$. Let $b_+(t)$ be the conditioned Brownian bridge and let $X_T(t) = b_+(t) - c_T(t)$ be the deviation of $b_+(t)$ away from $c_T(t)$, see Figure 1. Clearly $X_T(t) \geq 0$, $X_T(-T) = X_T(T) = 0$, and the path measure of the process is defined on $C([-T, T], \mathbb{R}) = C([-T, T])$, the space of continuous functions over the interval $[-T, T]$ equipped with the supremum norm. The issue is to understand the statistical properties of $X_T(t)$ for large $T$.

A well-studied special case is when $c_T(t)$ is replaced by the function zero. The Brownian bridge is then constrained to stay positive, a stochastic process known as Brownian excursion. In the limit of large $T$ it converges to the three-dimensional Bessel process. Time-dependent barriers, like the circle, seem to be hardly studied. An exception is the parabola $g_{T,2}(t) = T^2 - t^2$ for

Received August 2003; revised October 2004.

[1]Supported by Deutsche Forschungsgemeinschaft Research Project SP 181/17–1.

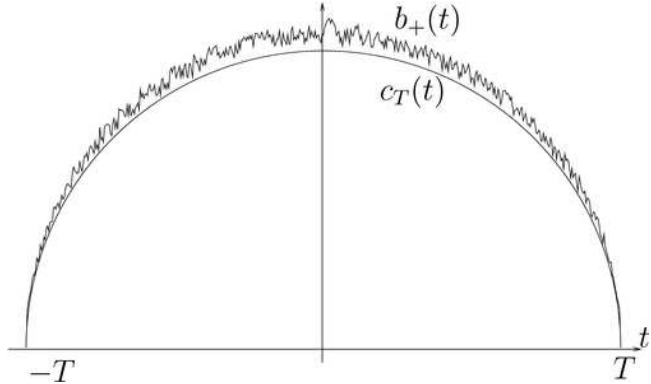

Fig. 1. *Brownian bridge $b_+(t)$ conditioned to lie above the semicircle $c_T(t)$.*

which some properties have been established [4, 5]; see below. In this paper we resolve the fluctuation problem for:

(i) the circle $c_T(t)$,
(ii) the family of parabolas $g_{T,\gamma} = T^\gamma(1 - (t/T)^2)$.

We also discuss briefly general shape functions of the form $g_T(t) = Tg(t/T)$.

Our problem arose rather indirectly in an attempt to understand a one-layer approximation to the multilayer polynuclear growth model; see [6]. There one has $N + 1$ independent copies of the Brownian bridge, denoted here as $b_j(t)$, $|t| \leq T$, $j = 0, -1, \ldots, -N$, such that $b_j(\pm T) = j$, and conditions them on nonintersection, with the subsequent limit $N \to \infty$. Of interest is the top line $b_0(t)$, $|t| \leq T$. Because of conditioning, typically $b_0(t)$ has a shape of a semicircle. Therefore the crude approximation consists in replacing all lower-lying Brownian motions, that is, $b_j(t)$ with $j = -1, -2, \ldots$, by the semicircle $c_T$. As we will prove, this approximation preserves the scaling behavior, in the sense that transverse fluctuations are of order $T^{1/3}$ and longitudinal correlations decay over a time span of order $T^{2/3}$. However, finer details are not accounted for. For example, in our problem $X_T(t)$, on the scale $T^{2/3}$, is exponentially mixing, whereas the covariance of top line $b_0(t)$ on the same scale has only power law decay [2, 10].

To state our main result we define the *stationary* diffusion process $\mathcal{A}(t)$ through the stochastic differential equation

$$(1.1) \qquad d\mathcal{A}(t) = a(\mathcal{A}(t))\,dt + db_t$$

with $b_t$ the standard Brownian motion and drift

$$(1.2) \qquad a(x) = \frac{\text{Ai}'(-\omega_1 + x)}{\text{Ai}(-\omega_1 + x)},$$

where $-\omega_1$ is the first zero of the Airy function Ai [1]. The relevant asymptotic is $a(x) = x^{-1}$ for $x \to 0^+$ and $a(x) = -\sqrt{x}$ for $x \to \infty$. Thus (1.1) admits a unique stationary measure which is given by

$$\frac{d}{dx}\mathbb{P}(\mathcal{A}(t) \leq x) = \frac{\text{Ai}(-\omega_1 + x)^2}{\text{Ai}'(-\omega_1)^2}\mathbb{1}_{[x>0]}. \tag{1.3}$$

$\mathcal{A}(t)$ has continuous sample paths and the small $x$ behavior of the drift implies that $\mathbb{P}(\mathcal{A}(t) > 0 \text{ for all } t) = 1$.

**Theorem 1.1.** *Let $b_+(t)$ be the Brownian bridge $b(t)$ conditioned on the set $\{b(t) \geq c_T(t) \text{ for all } t \in [-T,T]\}$ and let $X_T(t) = b_+(t) - c_T(t)$, $|t| \leq T$. The rescaled process close to the reference point $\tau T$ is defined through*

$$t \mapsto \mathcal{A}_T(t) = v_s X_T(\tau T + h_s^{-1}t), \tag{1.4}$$

*with $v_s = 2^{1/3}(1-\tau^2)^{-1/2}T^{1/3}$, $h_s = v_s^2$. Then*

$$\lim_{T \to \infty} \mathcal{A}_T = \mathcal{A}, \tag{1.5}$$

*in the sense of weak convergence of path measures on $C([-N,N])$, for any $N > 0$.*

For the polynuclear growth model, the same rescaling leads to the Airy process, which has a $t^{-2}$ decay of correlations as is known from the rather intricate explicit solution given in [2, 10]. This behavior should be seen in contrast to the exponential mixing of the diffusion process $\mathcal{A}(t)$.

To prove Theorem 1.1, we rely on the fact that some reasonably explicit expressions are available in case the semicircle is replaced by a parabola of the form

$$g_{T,\gamma}(t) = T^\gamma(1 - (t/T)^2). \tag{1.6}$$

**Theorem 1.2.** *Let $b_{+,\gamma}(t)$ be the Brownian bridge $b(t)$ conditioned on the set $\{b(t) \geq g_{T,\gamma}(t) \text{ for all } t \in [-T,T]\}$ and let $X_{T,\gamma}(t) = b_{+,\gamma}(t) - g_{T,\gamma}(t)$. The rescaled process is defined through*

$$t \mapsto \mathcal{A}_{T,\gamma}(t) = v_s X_{T,\gamma}(\tau T + h_s^{-1}t), \tag{1.7}$$

*with $v_s = T^{(\gamma-2)/3}4^{1/3}$, $h_s = v_s^2$. Then, for $\gamma > 1/2$,*

$$\lim_{T \to \infty} \mathcal{A}_{T,\gamma} = \mathcal{A}, \tag{1.8}$$

*in the sense of weak convergence of path measures on $C([-N,N])$, for any $N > 0$.*

The limit (1.7) has the, at first sight surprising, feature that the limit process $\mathcal{A}(t)$ does not depend on the scaling exponent $\gamma$. For $\gamma = 2$, that is, the standard parabola $g_{T,2}(t) = T^2 - t^2$, the fluctuations are of order 1, whereas for $\gamma > 2$ they actually decrease as $T \to \infty$. The condition $\gamma > 1/2$ reflects the fact that as $\gamma \to 1/2$ the time-scaling $T^{-2(\gamma-2)/3} \to T$. In other words, for $\gamma = 1/2$ the interior is correlated with the end-points and no stationary distribution is reached locally. For $\gamma < 1/2$, $g_{T,\gamma}(t)$ can be replaced by the function zero and the limit process is the Brownian excursion.

We outline the strategy to prove Theorem 1.1. Note that $X_T(t)$ is Markov, in the sense that upon conditioning on $X_T(t_0)$ the future and the past path measures are independent. Let us fix then the time window $[-N, N]$ for the rescaled process $\mathcal{A}_T(t)$.

(i) The first step is to show that the entrance/exit law, that is, the joint distribution of $(\mathcal{A}_T(-N), \mathcal{A}_T(N))$, is close to the corresponding entrance/exit law of the limit diffusion process $\mathcal{A}$. To achieve such a result the true shape function $c_T(t)$ is piecewise approximated by parabolas. Parabolas are chosen because for them reasonably explicit expressions for the transition probability are available.

(ii) For the interval $[-N, N]$ we use the limit entrance/exit law and use a suitably chosen parabola as conditioning shape, such that the resulting process is identical to $\mathcal{A}(t)$, $|t| \leq N$. Thus the claim of Theorem 1.1 follows from the fact that inside $[-N, N]$ the circle and the parabola differ at most by $\mathcal{O}(T^{-1/4})$.

Following this strategy, in Section 2 we consider the parabolic constraint and prove Theorem 1.2. In Section 3 we establish a result needed to control the joint entrance/exit law for the time window under consideration. With this input we prove Theorem 1.1 in Section 4. In Section 5 we discuss other shapes. The Appendix contains estimates on the transition probability for the conditioning parabolic constraint and some monotonicity results required in Section 4.

## 2. Parabolic constraint.

We plan to prove Theorem 1.2 and first state a result on the transition density for Brownian motion conditioned to remain below a parabola $-\frac{1}{2} g_{T,2}(t + T)$. This result was first obtained by Groeneboom; see (2.23) and (2.24) in [5]. In a different way it was derived by Salminen; see Proposition (3.9) of [8]. We were led to the explicit formula in Lemma 2.1 below from Frachebourg and Martin, page 330 of [4], where the references to [5, 8] are given. Since the result holds for an arbitrary diffusion coefficient, by Brownian motion scaling we can easily deduce the transition density for Brownian motion conditioned to remain above $g_{T,\gamma}(t)$. The result is reported in Lemma 2.1 below. The vertical and horizontal scaling depends only on the $g''_{T,\gamma}(t)$; therefore we define

$$(2.1) \qquad \kappa = -\frac{d^2 g_{T,\gamma}(t)}{dt^2} = 2T^{\gamma-2}.$$

Let $W(x_2, t_2 | x_1, t_1)$ be the transition probability density for Brownian motion $b^{x_1, t_1}(t)$ conditioned to start at $t_1$ from $g_{T,\gamma}(t_1) + x_1$ and ending at $t_2$ in $g_{T,\gamma}(t_2) + x_2$,

$$(2.2) \qquad W(x_2, t_2 | x_1, t_1) = \frac{d}{dx_2} \mathbb{P}(Y(t_2) \le x_2 | Y(t) \ge 0, t \in [t_1, t_2])$$

where $Y(t) = b^{x_1, t_1}(t) - g_{T,\gamma}(t)$.

LEMMA 2.1. *Let us define the vertical and horizontal scaling as*

$$(2.3) \qquad v_s = (2\kappa)^{1/3}, \qquad h_s = (2\kappa)^{2/3}.$$

*Then*

$$(2.4) \qquad W(x_2, t_2 | x_1, t_1) = \widehat{W}(x_2, t_2 | x_1, t_1) \exp(F(x_2, t_2 | x_1, t_1))$$

*with*

$$(2.5) \qquad \widehat{W}(x_2, t_2 | x_1, t_1) = \sum_{k \ge 1} v_s e^{-\omega_k (t_2 - t_1) h_s / 2} \frac{\mathrm{Ai}(v_s x_1 - \omega_k)}{\mathrm{Ai}'(-\omega_k)} \frac{\mathrm{Ai}(v_s x_2 - \omega_k)}{\mathrm{Ai}'(-\omega_k)}$$

*and*

$$(2.6) \qquad \begin{aligned} F(x_2, t_2 | x_1, t_1) &= x_1 g'_{T,\gamma}(t_1) - x_2 g'_{T,\gamma}(t_2) \\ &\quad - \frac{1}{6\kappa^{1/3}} [g'_{T,\gamma}(t_1)^3 - g'_{T,\gamma}(t_2)^3]. \end{aligned}$$

*Here* $-\omega_1, -\omega_2, \ldots$ *are the zeros of the Airy function,* $0 < \omega_1 < \omega_2 < \cdots$.

Let $X_{T,\gamma}(t)$ be the process of Theorem 1.2. Furthermore let $L$ be the backward generator of the diffusion process $\mathcal{A}(t)$,

$$(2.7) \qquad (L\varphi)(x) = \frac{1}{2} \frac{d^2 \varphi(x)}{dx^2} + a(x) \frac{d\varphi(x)}{dx}$$

as acting on smooth functions $\varphi$. $\mathcal{A}(t)$ has the invariant measure $\Omega(x)^2$ with

$$(2.8) \qquad \Omega(x) = \frac{\mathrm{Ai}(-\omega_1 + x)}{\mathrm{Ai}'(-\omega_1)}, \qquad x \ge 0, \qquad \int_{\mathbb{R}_+} \Omega(x)^2 \, dx = 1.$$

Through the ground-state transformation $Hf = -\Omega(L\Omega^{-1}\varphi)$ (see, e.g., Chapter V.16 of [9]), one obtains

$$(2.9) \qquad (H\varphi)(x) = -\frac{1}{2} \frac{d^2 \varphi(x)}{dx^2} + \frac{x}{2} \varphi(x) - E\varphi(x), \qquad x \ge 0.$$

$H$ is understood with Dirichlet boundary condition at $x = 0$ and $E = \frac{1}{2}\omega_1$ implies $H\Omega = 0$. Denote by $G(x, y; t)$ the integral kernel of $G_t = e^{-tH}$, that is,

$$(2.10) \qquad (e^{-tH}\varphi)(x) = \int_{\mathbb{R}_+} G(x, y; t)\varphi(y) \, dy.$$

We remark that $H$ has purely discrete spectrum. Its eigenvalues and eigenfunctions are given by

$$(2.11) \quad E_k = \tfrac{1}{2}\omega_{k+1}, \qquad \Omega_k(x) = \frac{\mathrm{Ai}(-\omega_{k+1} + x)}{\mathrm{Ai}'(-\omega_{k+1})}, \qquad x \geq 0, k = 0, 1, \ldots.$$

Note that we use the notation $\Omega \equiv \Omega_0$, since $\Omega_0$ reappears frequently throughout the paper.

Before proving Theorem 1.2 we explain how $\mathcal{A}$ is related to a conditioned Brownian motion.

PROPOSITION 2.2. *Let $Z(t)$, $|t| \leq N$, be Brownian motion conditioned to stay above $s(t) = \frac{1}{4}(N^2 - t^2)$ and such that the joint probability density of $(Z(-N), Z(N))$ is given by*

$$(2.12) \quad \rho_Z(\xi_1, -N; \xi_2, N) = \Omega(\xi_2)G(\xi_2, \xi_1; 2N)\Omega(\xi_1).$$

*Then $Z \overset{\mathcal{D}}{=} \mathcal{A} + s$ on $C([-N, N])$.*

PROOF. Denote $W(t) = Z(t) - s(t)$; then the transition density of $W(t)$ is

$$
\begin{aligned}
p(y, t | x, u) = \Bigg[ & \int_{\mathbb{R}_+^2} d\xi_1 \, d\xi_2 \, \rho_Z(\xi_1, -N; \xi_2, N) \\
& \times G(\xi_2, y; N - t)G(y, x; t - u) \\
& \times G(x, \xi_1; u + N)G(\xi_2, \xi_1; 2N)^{-1} \Bigg] \\
& \times \Bigg[ \int_{\mathbb{R}_+^2} d\xi_1 \, d\xi_2 \, \rho_Z(\xi_1, -N; \xi_2, N) \\
& \times G(\xi_2, x; N - u)G(x, \xi_1; u + N)G(\xi_2, \xi_1; 2N)^{-1} \Bigg]^{-1}
\end{aligned}
$$
(2.13)

for $x, y > 0$ and $-N < u < t < N$. But since $\rho_Z(\xi_1, -N; \xi_2, N) = \Omega(\xi_2) \times G(\xi_2, \xi_1; 2N)\Omega(\xi_1)$, it follows that

$$(2.14) \quad p(y, t | x, u) = (G_{N-t}\Omega)(y)G(y, x; t - u)/(G_{N-u}\Omega)(x).$$

Notice that $h(x) = (G_{N-t}\Omega)(x) = \Omega(x)$. Hence the process with transition probability density (2.14) is the Doob $h$-transform; see Section IV.39 of [7]. Thus it follows that the process $W(t)$ satisfies the SDE

$$(2.15) \quad dW(t) = \tilde{a}(W(t)) \, dt + db_t$$

with the drift $\tilde{a}(x) = \partial \ln h(x)/\partial x$ being equal to (1.2) and $b_t$ standard Brownian motion. Therefore $W(t)$ and $\mathcal{A}(t)$ satisfy the same SDE and, since they have the same distribution at $t = -N$, $W(t) \overset{\mathcal{D}}{=} \mathcal{A}(t)$. □

We now prove Theorem 1.2 for the case of the parabolic constraint $g_{T,\gamma}$. The strategy consists in first controlling the joint density of $(\mathcal{A}_{T,\gamma}(-N), \mathcal{A}_{T,\gamma}(N))$, and then using the Markov property of Brownian motion together with Proposition 2.2 to determine the limit process of $\mathcal{A}_{T,\gamma}$. This strategy will be also the basis of the proof of Theorem 1.1.

PROOF OF THEOREM 1.2. Consider the rescaled process $\mathcal{A}_{T,\gamma} = v_s \times X_{T,\gamma}(\tau T + h_s^{-1}t)$, $|t| \leq N$, with $v_s = T^{(\gamma-2)/3}4^{1/3}$ and $h_s = v_s^2$. The joint density of $(\mathcal{A}_{T,\gamma}(-N), \mathcal{A}_{T,\gamma}(N))$ is given by

$$
\begin{aligned}
\rho_T(\xi_1, -N; \xi_2, N) = \lim_{\varepsilon \to 0} &G(\varepsilon, \xi_2; T(1-\tau)h_s - N)G(\xi_2, \xi_1; 2N) \\
&\times G(\xi_1, \varepsilon; T(1+\tau)h_s - N)/G(\varepsilon, \varepsilon; 2Th_s).
\end{aligned}
\tag{2.16}
$$

Since $\gamma > 1/2$, $Th_s \sim T^{(2\gamma-1)/3} \to \infty$ when $T \to \infty$. Using the estimate from Lemma A.1, we have, for some constant $a > 0$,

$$
G(\varepsilon, \varepsilon; 2Th_s) = \varepsilon^2(1 + \mathcal{O}(e^{-aTh_s}))
\tag{2.17}
$$

and

$$
\begin{aligned}
&G(\varepsilon, \xi_2; T(1-\tau)h_s - N) \\
&= \varepsilon[\Omega(\xi_2) + \mathcal{O}(\min\{\xi_2 e^{-aTh_s}, e^{-a\xi_2(Th_s)^{1/3}}\})].
\end{aligned}
\tag{2.18}
$$

Therefore

$$
\begin{aligned}
\lim_{T \to \infty} \rho_T(\xi_1, -N; \xi_2, N) &= \Omega(\xi_2)G(\xi_2, \xi_1; 2N)\Omega(\xi_1) \\
&\equiv \rho_{\mathcal{A}}(\xi_1, -N; \xi_2, N).
\end{aligned}
\tag{2.19}
$$

For any bounded, continuous function $f$ on $C([-N, N])$,

$$
\begin{aligned}
\mathbb{E}_{\mathcal{A}_{T,\gamma}}(f) &= \int_{\mathbb{R}_+^2} d\xi_1 \, d\xi_2 \, \rho_T(\xi_1, -N; \xi_2, N) \\
&\quad \times \mathbb{E}_{\mathcal{A}_{T,\gamma}}(f|\mathcal{A}_{T,\gamma}(-N) = \xi_1, \mathcal{A}_{T,\gamma}(N) = \xi_2) \\
&= \int_{\mathbb{R}_+^2} d\xi_1 \, d\xi_2 \, \rho_{\mathcal{A}}(\xi_1, -N; \xi_2, N) \\
&\quad \times \mathbb{E}_{\mathcal{A}_{T,\gamma}}(f|\mathcal{A}_{T,\gamma}(-N) = \xi_1, \mathcal{A}_{T,\gamma}(N) = \xi_2) \\
&\quad + R_1(T, N, f),
\end{aligned}
\tag{2.20}
$$

with $R_1(T, N, f)$ bounded by

$$
\begin{aligned}
|R_1(T, N, f)| &\leq \|f\|_\infty \int_{\mathbb{R}_+^2} d\xi_1 \, d\xi_2 \\
&\quad \times |\rho_T(\xi_1, -N; \xi_2, N) - \rho_{\mathcal{A}}(\xi_1, -N; \xi_2, N)|,
\end{aligned}
\tag{2.21}
$$

which converges to zero as $T \to \infty$, because $\rho_T$ converges pointwise to $\rho_{\mathcal{A}}$ and $\rho_T$, $\rho_{\mathcal{A}}$ are densities with total mass 1 (Scheffé's theorem; see, e.g., Appendix of [3]). Finally, Proposition 2.2 implies that the nonvanishing term in (2.20) is $\mathbb{E}_{\mathcal{A}}(f)$.  $\square$

**3. Joint entrance and exit law.** In a piecewise parabolic approximation of the semicircle, or more generally of a *concave* function, there are points with discontinuities in the slope. In order to control the subleading terms we take a continuous, piecewise parabolic shape such that the *derivative has negative jumps* at its discontinuity points. We call these points *ridges*.

More precisely, let us consider a Brownian bridge $b_s(t)$ conditioned to remain above a continuous, concave, piecewise parabolic function $s(t)$, starting from $s(t_{\text{in}}) + x_{\text{in}}$ at time $t_{\text{in}}$ and ending at $s(t_{\text{fin}}) + x_{\text{fin}}$ at time $t_{\text{fin}}$, $t_{\text{in}} < t_{\text{fin}}$, where

$$(3.1) \qquad s(t) = a_j + b_j t - \tfrac{1}{2} c_j t^2 \qquad \text{for } t \in [u_{j-1}, u_j]$$

with $c_j > 0$, $u_0 = t_{\text{in}}$ and $u_{M-1} = t_{\text{fin}}$. We want to study the process close to $t = \tilde{t}$, with $\tilde{t}$ very far away from the contact times $u_j$, say $u_{K-1} \ll \tilde{t} \ll u_K$. Define

$$(3.2) \qquad v_s = (-2s''(\tilde{t}))^{1/3}, \qquad h_s = v_s^2,$$

the times $t_j = u_j$ for $j = 0, \dots, K-1$, $t_j = u_{j-1}$ for $j = K+1, \dots, M$, and

$$(3.3) \qquad t_K \equiv t_- = \tilde{t} - N h_s^{-1}, \qquad t_{K+1} \equiv t_+ = \tilde{t} + N h_s^{-1}.$$

Denote

$$(3.4) \qquad \nu(t_j) = s'(t_j^-) - s'(t_j^+) \geq 0,$$

in particular, $\nu(t_K) = \nu(t_{K+1}) = 0$, and

$$(3.5) \quad v_j = (-2s''((t_j + t_{j-1})/2))^{1/3}, \qquad h_j = v_j^2, \qquad \Gamma_j = \tfrac{1}{2} h_j(t_j - t_{j-1}).$$

Finally, let $\bar{\Gamma} = \min_{j \neq K+1} \Gamma_j$, assume that $\bar{\Gamma} \to \infty$ as $T \to \infty$, and that $v_1 x_{\text{in}} \leq \Gamma_1^{1/3}$, $v_M x_{\text{fin}} \leq \Gamma_M^{1/3}$.

LEMMA 3.1. *Fix $N > 0$ and denote $t \mapsto X_T(t) = b_s(t) - s(t)$. Then the joint density of $(X_T(t_-), X_T(t_+))$ with $t_- = t_K - N h_s^{-1}$ and $t_+ = t_K + N h_s^{-1}$ is given by*

$$\rho_T(x, t_-; y, t_+) = \frac{d}{dx}\frac{d}{dy}\mathbb{P}(X_T(t_-) \leq x, X_T(t_+) \leq y)$$

$$(3.6) \qquad = v_s^2 \Omega(v_s x) G(v_s x, v_s y; 2N) \Omega(v_s y)(1 + \mathcal{O}(e^{-2a\bar{\Gamma}^{1/3}}))$$

$$\qquad + E_T(x, t_-; y, t_+),$$

*for some constant $a > 0$ and where the error term $E_T$ converges pointwise to 0 and its total mass is bounded by*

$$(3.7) \qquad \int_{\mathbb{R}_+^2} dx\, dy\, |E_T(x, t_-; y, t_+)| = \mathcal{O}(e^{-a\Gamma^{1/3}}).$$

PROOF. Let us denote by $z_i$ the position of the Brownian bridge above $s(t_i)$ for $i = 0, \ldots, M$. Then the density (3.6) is given by

$$(3.8) \qquad \rho_T(x, t_-; y, t_+) = \frac{\int_{\mathbb{R}_+^{M-3}} \prod_{i \in J} dz_i \prod_{j=1}^{M} W(z_j, t_j | z_{j-1}, t_{j-1})}{\int_{\mathbb{R}_+^{M-1}} \prod_{i=1}^{M-1} dz_i \prod_{j=1}^{M} W(z_j, t_j | z_{j-1}, t_{j-1})}$$

with $J = \{1, \ldots, K-1, K+2, \ldots, M-1\}$. Explicitly

$$(3.9) \qquad \begin{aligned} W(z_j, t_j | z_{j-1}, t_{j-1}) &= \widehat{W}(z_j, t_j | z_{j-1}, t_{j-1}) \\ &\quad \times \exp[z_{j-1} s'(t_{j-1}^+) - z_j s'(t_j^-)] q(t_j, t_{j-1}) \end{aligned}$$

with $q$ a function independent of $z_j, z_{j-1}$. When (3.9) is substituted in (3.8), the product of the $q$'s simplifies. Moreover, each $\widehat{W}$ contains a prefactor $v_j e^{-\omega_1 \Gamma_j}$; see (A.5). Thus $W(z_j, t_j | z_{j-1}, t_{j-1})$ in (3.8) can be replaced by

$$(3.10) \qquad v_j^{-1} e^{\omega_1 \Gamma_j} \widehat{W}(z_j, t_j | z_{j-1}, t_{j-1}) \exp[z_{j-1} s'(t_{j-1}^+) - z_j s'(t_j^-)]$$

and in addition $s'(t_0^+)$ and $s'(t_{M+1}^-)$ can be replaced by zero.

Let us first analyze the denominator of (3.8). It can be written

$$(3.11) \quad \int_{\mathbb{R}_+^{M-1}} \prod_{i=1}^{M-1} (dz_i\, e^{-z_i \nu(t_i)}) \prod_{j=1}^{M} (\Omega(v_j z_j)\Omega(v_j z_{j-1}) + R_{\Gamma_j}(v_j z_j, v_j z_{j-1})),$$

where $R_{\Gamma_j}$ is the one in Lemma A.1. Denote

$$(3.12) \quad Q = \Omega(v_1 z_0)\Omega(v_M z_M) \prod_{i \in J} \int_{\mathbb{R}_+} dz_i\, e^{-z_i \nu(t_i)} \Omega(v_i z_i)\Omega(v_{i+1} z_i);$$

then the expansion of (3.11) has the leading term

$$(3.13) \qquad Q \prod_{i=K}^{K+1} \int_{\mathbb{R}_+} dz_i\, \Omega(v_s z_i)\Omega(v_s z_i) = Q v_s^{-2}$$

plus $2^M - 1$ terms containing one or more factors of $R$'s. The conditions $v_1 x_{\text{in}} \le \Gamma_1^{1/3}$ and $v_M x_{\text{fin}} \le \Gamma_M^{1/3}$ imply the bounds

$$(3.14) \qquad \begin{aligned} |R_{\Gamma_1}^0(v_1 x_{\text{in}})| &\le \Omega(v_1 x_{\text{in}})\mathcal{O}(e^{-a\Gamma_1^{2/3}/2}), \\ |R_{\Gamma_M}^0(v_M x_{\text{fin}})| &\le \Omega(v_M x_{\text{fin}})\mathcal{O}(e^{-a\Gamma_M^{2/3}/2}). \end{aligned}$$

Using Lemma A.2, we can replace each $R_{\bar{\Gamma}}^0$ by $\Omega$ in the integration variables up to a multiplicative factor $\mathcal{O}(e^{-a\Gamma^{1/3}})$. Summing up all these contributions, the denominator is given by

$$(3.15) \qquad \text{denominator of (3.8)} = \alpha Q v_s^{-2}(1 + \mathcal{O}(e^{-2a\bar{\Gamma}^{1/3}}))$$

where $\alpha \neq 0$ is a constant coming from the replacements described before (3.11).

The numerator is obtained similarly, but the variables $x$ and $y$ are not integrated out, with the result

$$(3.16) \quad \text{numerator of (3.8)} = \alpha Q \Omega(v_s y) G(v_s y, v_s x; 2N) \Omega(v_s x) + E_1(x, y)$$

where the first is the term with no factor of $R$ and $E_1(x, y)$ is the error term, which is bounded by

$$
\begin{aligned}
(3.17) \quad |E_1(x,y)| \leq{}& \alpha Q G(v_s y, v_s x; 2N) \\
&\times [\Omega(v_s x) R_{\bar{\Gamma}}^0(v_s y) + \Omega(v_s y) R_{\bar{\Gamma}}^0(v_s x) + R_{\bar{\Gamma}}^0(v_s x) R_{\bar{\Gamma}}^0(v_s y)] \\
&\times (1 + \mathcal{O}(e^{-2a\bar{\Gamma}^{1/3}}))
\end{aligned}
$$

with $R_{\bar{\Gamma}}^0$ given in (A.6). From (3.15) and (3.16) it follows that

$$
\begin{aligned}
(3.18) \quad \rho_T(x, t_-; y, t_+) ={}& v_s^2 \Omega(v_s y) G(v_s y, v_s x; 2N) \Omega(v_s x)(1 + \mathcal{O}(e^{-2a\bar{\Gamma}^{1/3}})) \\
&+ E_2(x, y)
\end{aligned}
$$

with $E_2(x, y) = E_1(x, y)/Q v_s^{-2}(1 + \mathcal{O}(e^{-2a\bar{\Gamma}^{1/3}}))$.

The expression of $R_{\bar{\Gamma}}^0$ implies that $R_{\bar{\Gamma}}^0(y) \leq e^{-a\bar{\Gamma}^{1/3}}$, converges pointwise to 0, and decays exponentially in $y$ for large $y$. On the other hand, $G(y, x; 2N)$ is uniformly bounded in $x$ and $y$ for any $N > 0$. Therefore

$$(3.19) \qquad \int_{\mathbb{R}_+^2} dx\, dy\, |E_2(x, y)| \leq \mathcal{O}(e^{-a\bar{\Gamma}^{1/3}}). \qquad \qquad \square$$

## 4. Proof of Theorem 1.1.

In order to prove the theorem we first control the entrance/exit law for the interval $[\tau T - N h_s^{-1}, \tau T + N h_s^{-1}]$, for which we use Lemma 3.1. Therefore one has to find a lower and an upper approximation satisfying its hypotheses.

4.1. *Upper and lower approximating shapes for $t = -\tau T$*. The piecewise parabolic approximations $s_\pm$ are constructed with the parabolas

$$(4.1) \qquad f_i(t) = a_i + b_i t - \tfrac{1}{2} c_i t^2 \qquad \text{for } t \in [u_{i-1}, u_i]$$

for $-T = u_0 < u_1 < \cdots < u_{n-1} < u_n = 0$, where the choice of the $u_j$'s is discussed below. We set $s_\pm(t) = s_\pm(-t)$ for $t \in [0, T]$ (although this is not required for the result). Since we want to apply Proposition 2.2, we also determine $v_j = (2c_j)^{1/3}$ and $\Gamma_j = \frac{1}{2}(u_j - u_{j-1})h_j$ with $h_j = v_j^2$. In case $\tau = 0$, we set $b_j = 0$.

4.1.1. *Upper approximation, $\tau = 0$.* This is the easiest case and one needs only a single parabola, that is, $n = 1$,

$$(4.2) \qquad f_1(t) = T - \tfrac{1}{2}T^{-1}t^2.$$

$s_+(t) = f_1(t) \geq c_T(t)$ for all $t$. Since $u_0 = -T$, $u_1 = 0$,

$$(4.3) \qquad v_1 = 2^{1/3}T^{-1/3}, \qquad \Gamma_1 = 2^{-1/3}T^{1/3}.$$

4.1.2. *Lower approximation, $\tau = 0$.* In this case one needs $n = 2$. We define $u_1 = -T^{3/4}$. The parabola from $(-T, 0)$ to $(u_1, c_T(u_1))$ is given by (4.1) with $a_1 = T(1 - T^{-1/2})^{-1/2}$ and $c_1 = 2T^{-1} + \mathcal{O}(T^{-3/2})$. The parabola from $(u_1, c_T(u_1))$ to $(0, T)$ has $a_2 = T$ and $c_2 = T^{-1} + \mathcal{O}(T^{-3/2})$. Then for $t \in [-T, T]$, $s_-(t) \leq c_T(t)$, with

$$(4.4) \qquad v_1 = 2^{2/3}T^{-1/3} + \mathcal{O}(T^{-5/6}), \qquad \Gamma_1 = 2^{1/3}T^{1/3} + \mathcal{O}(T^{1/12})$$

and

$$(4.5) \qquad v_2 = 2^{1/3}T^{-1/3} + \mathcal{O}(T^{-5/6}), \qquad \Gamma_2 = 2^{-1/3}T^{1/12} + \mathcal{O}(T^{-5/12}).$$

$s_-(t)$ has a ridge at $\pm u_1$.

4.1.3. *Upper approximation, $\tau < 0$.* In this case the construction requires $n = 3$. For convenience we define $\lambda_\tau = 1 - \tau^2$ and $\beta = -\tau > 0$. Let $u_1 = -\tau T$ and let the parabola $f_1(t)$ be defined by

$$(4.6) \qquad \begin{aligned} f_1(t) = f_2(t) &= c_T(\tau T) + c_T'(\tau T)(t - \tau T) \\ &\quad + \tfrac{1}{2}c_T''(\tau T)(1 - T^{-1/4})(t - \tau T)^2. \end{aligned}$$

We define $u^*$ to be the first intersection time after $u_1$ of $f_2(t)$ with $c_T(t)$. We estimate $u^* = -\beta T + \lambda_\tau \beta^{-1}T^{3/4} + \mathcal{O}(T^{1/2})$. Let

$$(4.7) \qquad f^*(t) = a^* - \tfrac{1}{2}c^*t^2$$

be the parabola which passes through $(u_1, c_T(t_1))$ and $(u^*, c_T(u^*))$. Some computations lead to $c^* = \lambda_\tau^{-1/2}T^{-1} + \mathcal{O}(T^{-5/4})$. Since $f^*(t) \leq c_T(t)$ for $t \in [u_1, u^*]$ and $f_2(t) \geq c_T(t)$ for $t \in [u_1, u^*]$, there is a time $u_2 \in (u_1, u^*)$ such that $f_2'(u_2) = f^{*'}(u_2)$. We obtain $u_2 = -\beta T + \tfrac{1}{2}\lambda_\tau \beta^{-1}T^{3/4} + \mathcal{O}(T^{1/2})$. Finally one has to define the third piece of parabola. Since $f^*(t) \geq c_T(t)$ for $t \geq u^*$, and $f_2(t) \geq c_T(t)$ for $t \in [u_1, u^*]$, we define $f_3(t)$ by

$$(4.8) \qquad f_3(t) = f^*(t) + (f_2(u_2) - f^*(u_2)).$$

This construction satisfies $s_+(t) \geq c_T(t)$ for $t \in [-T, T]$, has a ridge at $t = 0$, and the second derivative is discontinuous at $t = \pm u_2$. Moreover one has

$$
\begin{aligned}
(4.9) \quad & v_1 = 2^{1/3} \lambda_\tau^{-1/2} T^{-1/3} + \mathcal{O}(T^{-7/12}), \\
& \Gamma_1 = 2^{-1/3}(1 - |\tau|)^{-1} T^{1/3} + \mathcal{O}(T^{1/12}), \\
& v_2 = v_1, \\
& \Gamma_2 = 2^{-4/3}|\tau|^{-1} T^{1/12} + \mathcal{O}(T^{-1/6}), \\
& v_3 = 2^{1/3} \lambda_\tau^{-1/6} T^{-1/3} + \mathcal{O}(T^{-7/12}), \\
& \Gamma_3 = 2^{-1/3}|\tau| \lambda_\tau^{-1/3} T^{1/3} + \mathcal{O}(T^{1/12}).
\end{aligned}
$$

4.1.4. *Lower approximation*, $\tau < 0$.   In this case the construction requires $n = 4$. Also here let $\beta = -\tau$ and $\lambda_\tau = 1 - \tau^2$. We define $u_2 = -\tau T$ and the parabola $f_2(t)$ by

$$
\begin{aligned}
(4.10) \quad f_2(t) = f_3(t) = {} & c_T(\tau T) + c_T'(\tau T)(t - \tau T) \\
& + \tfrac{1}{2} c_T''(\tau T)(1 + T^{-1/4})(t - \tau T)^2.
\end{aligned}
$$

$f_2(t)$ has an intersection with $c_T(t)$ for some time $t < u_2$, which we define to be $u_1$, and remains below $c_T(t)$ for $t \in [u_2, 0]$. Some computations lead to $u_1 = -\beta T - \lambda_\tau \beta^{-1} T^{3/4} + \mathcal{O}(T^{1/2})$. Moreover let

$$
(4.11) \qquad f_1(t) = a_1 - \tfrac{1}{2} c_1 t^2
$$

be the parabola passing through $(-T, 0)$ and $(u_1, c_T(u_1))$. It has $c_1 = 2\lambda_\tau^{-1/2} \times T^{-1} + \mathcal{O}(T^{-5/4})$. Finally we define $u_3 = -\beta T(1 - T^{-1/4})$ and

$$
(4.12) \qquad f_4(t) = a_4 - \tfrac{1}{2} c_4 t^2
$$

such that $f_4(u_3) = f_3(u_3)$ and $f_4'(u_3) = f_3'(u_3)$. We obtain $c_4 = \lambda_\tau^{-1/2} T^{-1} + \mathcal{O}(T^{-5/4})$.

This construction satisfies $s_-(t) \leq c_T(t)$ for $t \in [-T, T]$, has a ridge at $t = 0$ and at $t = \pm u_1$, and the second derivative is discontinuous at $t = \pm u_3$. Moreover one has

$$
\begin{aligned}
& v_1 = 2^{1/3} \lambda_\tau^{-1/6} T^{-1/3} + \mathcal{O}(T^{-7/12}), \\
& \Gamma_1 = 2^{-1/3}(1 - |\tau|) \lambda_\tau^{-1/2} T^{1/3} + \mathcal{O}(T^{1/12}), \\
& v_2 = 2^{1/3} \lambda_\tau^{-1/2} T^{-1/3} + \mathcal{O}(T^{-7/12}), \\
& \Gamma_2 = 2^{-1/3}|\tau|^{-1} T^{1/12} + \mathcal{O}(T^{-1/6}), \\
& v_3 = v_2, \\
(4.13) \quad & \Gamma_3 = 2^{-1/3}|\tau| \lambda_\tau^{-1} T^{1/12} + \mathcal{O}(T^{-1/6}),
\end{aligned}
$$

$$v_4 = 2^{1/3}\lambda_\tau^{-1/6}T^{-1/3} + \mathcal{O}(T^{-7/12}),$$

$$\Gamma_4 = 2^{-1/3}|\tau|\lambda_\tau^{-1/3}T^{1/3} + \mathcal{O}(T^{1/12}).$$

4.2. *Joint densities.* We compute now the joint entrance/exit law for the process of Theorem 1.1.

Let $b_\pm(t)$ be the Brownian bridge from $(s_\pm(-T), -T)$ to $(s_\pm(T), T)$ conditioned to stay above $s_\pm$. The processes we actually want to study are

$$(4.14) \qquad \mathcal{A}_{T,\pm}(t) = v_c[b_\pm(\tau T + h_c^{-1}t) - c_T(\tau T + h_c^{-1}t)]$$

and Proposition 2.2 is concerned with the processes

$$(4.15) \qquad Y_{T,\pm}(t) = v_{s_\pm}[b_\pm(\tau T + h_{s_\pm}^{-1}t) - s_\pm(\tau T + h_{s_\pm}^{-1}t)].$$

Let us denote $\lambda_{T,\pm} = v_{s_\pm}/v_c$, and

$$(4.16) \qquad g_{T,\pm}(t) = v_c[s_\pm(\tau T + h_c^{-1}t) - c_T(\tau T + h_c^{-1}t)].$$

Then

$$(4.17) \qquad \mathcal{A}_{T,\pm}(t) = \lambda_{T,\pm}^{-1}Y_{T,\pm}(\lambda_{T,\pm}^2 t) + g_{T,\pm}(t).$$

We compute $\lambda_{T,\pm}$ and bound $g_{T,\pm}(t)$ for $t \in [-N, N]$ with the result:

(a) Case $\tau = 0$,

$$(4.18) \qquad \begin{aligned} &\lambda_{T,+} = 1, &&g_{T,+}(t) = \mathcal{O}(N^4 T^{-2/3}), \\ &\lambda_{T,-} = 1 + \mathcal{O}(T^{-1/2}), &&g_{T,-}(t) = \mathcal{O}(N^2 T^{-1/2}). \end{aligned}$$

(b) Case $\tau < 0$,

$$(4.19) \qquad \lambda_{T,\pm} = 1 + \mathcal{O}(T^{-1/4}), \qquad g_{T,\pm}(t) = \mathcal{O}(N^2 T^{-1/4}).$$

LEMMA 4.1. *Let* $\rho_{T,c_T}(\xi_1, -N; \xi_2, N)$ *be the joint probability density of* $(\mathcal{A}_T(-N), \mathcal{A}_T(N))$, *where* $\mathcal{A}_T$ *is defined in* (1.4). *Then*

$$(4.20) \qquad \lim_{T\to\infty} \rho_{T,c_T}(\xi_1, -N; \xi_2, N) = \rho_{\mathcal{A}}(\xi_1, -N; \xi_2, N)$$

*with*

$$(4.21) \qquad \rho_{\mathcal{A}}(\xi_1, -N; \xi_2, N) \equiv \Omega(\xi_2)G(\xi_2, \xi_1; 2N)\Omega(\xi_1).$$

PROOF. Let $\rho_{T,\pm}(\xi_1, -N; \xi_2, N)$ be the joint probability density of $(\mathcal{A}_{T,\pm}(-N), \mathcal{A}_{T,\pm}(N))$. Then, since $\lambda_{T,\pm} \to 1$ and $g_{T,\pm}(t) \to 0$ as $T \to \infty$,

$$(4.22) \qquad \begin{aligned} \lim_{T\to\infty} \rho_{T,\pm}(\xi_1, -N; \xi_2, N) &= \rho_{\mathcal{A}}(\xi_1, -N; \xi_2, N) \\ &\equiv \Omega(\xi_2)G(\xi_2, \xi_1; 2N)\Omega(\xi_1). \end{aligned}$$

Denote $F_{T,*}(\xi_1, -N; \xi_2, N) = \int_{x_i \leq \xi_i} dx_1\, dx_2\, \rho_{T,*}(x_1, -N; x_2, N)$, where $* = \{+, -, c_T\}$. From the monotonicity properties of Propositions A.5 and A.6 it follows that

$$(4.23) \quad F_{T,+}(\xi_1, -N; \xi_2, N) \leq F_{T,c_T}(\xi_1, -N; \xi_2, N) \leq F_{T,-}(\xi_1, -N; \xi_2, N).$$

Taking the limit $T \to \infty$ in (4.23) and using (4.22) we obtain

$$(4.24)\quad \begin{aligned} \lim_{T \to \infty} F_{T,c_T}(\xi_1, -N; \xi_2, N) &= F_{\mathcal{A}}(\xi_1, -N; \xi_2, N)\\ &\equiv \int_{x_i \leq \xi_i} dx_1\, dx_2\, \rho_{\mathcal{A}}(\xi_1, -N; \xi_2, N),\end{aligned}$$

thus also

$$(4.25) \qquad \lim_{T \to \infty} \rho_{T,c_T}(\xi_1, -N; \xi_2, N) = \rho_{\mathcal{A}}(\xi_1, -N; \xi_2, N). \qquad \square$$

Finally we are in position to prove our main theorem on the circular constraint.

PROOF OF THEOREM 1.1. The process we have to analyze is

$$(4.26) \qquad \mathcal{A}_T(t) = v_s X_T(\tau T + h_s^{-1} t)$$

where $X_T(t)$ is defined in Theorem 1.1. We have to prove that $\mathcal{A}_T \xrightarrow{\mathcal{D}} \mathcal{A}$ on $C([-N, N])$ in the limit $T \to \infty$, which is done through

$$(4.27) \qquad \mathcal{A}_T + \tilde{c}_T \xrightarrow{\mathcal{D}} \mathcal{A} + s$$

where $s$ is a fixed parabola and $\tilde{c}_T$ is a (nonrandom) function satisfying

$$(4.28) \qquad \lim_{T \to \infty} \sup_{t \in [-N, N]} |\tilde{c}_T(t) - s(t)| = 0.$$

Then (4.27) implies

$$(4.29) \qquad \mathcal{A}_T + \tilde{c}_T - s \xrightarrow{\mathcal{D}} \mathcal{A},$$

since the mapping $x \mapsto x - s$ is continuous. Finally (4.28) combined with (4.29) implies that $\mathcal{A}_T \xrightarrow{\mathcal{D}} \mathcal{A}$ as $T \to \infty$.

Now, let us prove (4.27). Define $L_T(t) = \alpha_T(t - \tau T) + \beta_T$ to be the line intersecting the circle $c_T$ at times $t = \tau T \pm h_s^{-1} N$. Moreover, let

$$(4.30) \qquad \tilde{c}_T(t) = v_s(c_T(\tau T + h_s^{-1} t) - L_T(\tau T + h_s^{-1} t))$$

and

$$(4.31) \qquad s(t) = \tfrac{1}{4}(N^2 - t^2).$$

A simple calculation shows that $\tilde{c}_T(t) = s(t) + \mathcal{O}(N^3 T^{-1/3})$, $t \in [-N, N]$.

We now consider the process $Y_T = \mathcal{A}_T + \tilde{c}_T$. Let $f$ be a bounded, continuous function on $C([-N, N])$. Using the Markov property,

$$
\begin{aligned}
\mathbb{E}_{Y_T}(f) &= \int_{\mathbb{R}^2_+} d\xi_1 \, d\xi_2 \, \rho_{T,c_T}(\xi_1, -N; \xi_2, N) \\
&\quad \times \mathbb{E}_{Y_T}(f | Y_T(-N) = \xi_1, Y_T(N) = \xi_2) \\
&= \int_{\mathbb{R}^2_+} d\xi_1 \, d\xi_2 \, \rho_{\mathcal{A}}(\xi_1, -N; \xi_2, N) \\
&\quad \times \mathbb{E}_{Y_T}(f | Y_T(-N) = \xi_1, Y_T(N) = \xi_2) \\
&\quad + R_1(T, N, f),
\end{aligned}
\tag{4.32}
$$

where the remainder term $R_1(T, N, f)$ can be bounded by

$$
\begin{aligned}
|R_1(T, N, f)| &\leq \|f\|_\infty \int_{\mathbb{R}^2_+} d\xi_1 \, d\xi_2 \\
&\quad \times |\rho_{T,c_T}(\xi_1, -N; \xi_2, N) - \rho_{\mathcal{A}}(\xi_1, -N; \xi_2, N)|
\end{aligned}
\tag{4.33}
$$

which converges to zero as $T \to \infty$, because $\rho_{T,c_T}$ converges by Lemma 4.1 pointwise to $\rho_{\mathcal{A}}$, and $\rho_{T,c_T}, \rho_{\mathcal{A}}$ are densities with total mass 1 (Scheffé's theorem; see, e.g., Appendix of [3]).

Let $Z(t)$ be the process defined in Proposition 2.2 with joint density of $(Z(-N), Z(N))$ given by $\rho_Z(\xi_1, -N; \xi_2, N) = \rho_{\mathcal{A}}(\xi_1, -N; \xi_2, N)$. For any realization $\omega$ of $Z$, define $\chi_{\tilde{c}_T}(\omega) = 1$ if $\omega(t) \geq \tilde{c}_T(t)$ for all $t \in [-N, N]$ and $\chi_{\tilde{c}_T}(\omega) = 0$ otherwise. Then the leading term of (4.32) is

$$
\mathbb{E}_Z(f \chi_{\tilde{c}_T}) / \mathbb{E}_Z(\chi_{\tilde{c}_T}),
\tag{4.34}
$$

and we have to show that it converges to $\mathbb{E}_Z(f \chi_s) / \mathbb{E}_Z(\chi_s)$ as $T \to \infty$. Notice that the reference measure does not depend on $T$; the only $T$-dependent quantity is $\tilde{c}_T$. It is easy to see that

$$
\frac{\mathbb{E}_Z(f \chi_s)}{\mathbb{E}_Z(\chi_s)} = \frac{\mathbb{E}_Z(f \chi_{\tilde{c}_T})}{\mathbb{E}_Z(\chi_{\tilde{c}_T})} + R_2(s, \tilde{c}_T, f)
\tag{4.35}
$$

with

$$
\begin{aligned}
R_2(s, \tilde{c}_T, f) &= \frac{\mathbb{E}_Z(f \chi_s (1 - \chi_{\tilde{c}_T}))}{\mathbb{E}_Z(\chi_s)} - \frac{\mathbb{E}_Z(f(1 - \chi_s) \chi_{\tilde{c}_T})}{\mathbb{E}_Z(\chi_s)} \\
&\quad + \frac{\mathbb{E}_Z(f \chi_{\tilde{c}_T})}{\mathbb{E}_Z(\chi_{\tilde{c}_T})} \left( \frac{\mathbb{E}_Z(\chi_{\tilde{c}_T}) - \mathbb{E}_Z(\chi_s)}{\mathbb{E}_Z(\chi_s)} \right).
\end{aligned}
\tag{4.36}
$$

Equation (4.36) can be bounded as

$$
\begin{aligned}
|R_2(s, \tilde{c}_T, f)| &\leq \frac{2\|f\|_\infty}{\mathbb{E}_Z(\chi_s)} (\mathbb{E}_Z(\chi_{\tilde{c}_T}(1 - \chi_s)) + \mathbb{E}_Z(\chi_s(1 - \chi_{\tilde{c}_T}))) \\
&= \frac{2\|f\|_\infty}{\mathbb{E}_Z(\chi_s)} \mathbb{P}_Z(\{\chi_s \neq \chi_{\tilde{c}_T}\}).
\end{aligned}
\tag{4.37}
$$

Let $B_T = \{\omega | \chi_s(\omega) \neq \chi_{\tilde{c}_T}(\omega)\}$; then $\mathbb{P}_Z(\{\chi_s \neq \chi_{\tilde{c}_T}\}) = \mathbb{P}_Z(B_T)$. Let $\varepsilon_T = \sup_{t \in [-N,N]} |\tilde{c}_T(t) - s(t)| = \mathcal{O}(?)$ then $B_T \subset D_T = \{\omega | \chi_{s-\varepsilon_T}(\omega) \neq \chi_{s+\varepsilon_T}(\omega)\}$. In the limit $T \to \infty$, $\omega \in B_T$ if $\omega$ touches without crossing the parabola $s$. Such paths have probability zero, therefore $\lim_{T \to \infty} \mathbb{E}_{Y_T}(f) = \mathbb{E}(f\chi_s)/\mathbb{E}(\chi_s)$.

We have proved that $Y_T = \mathcal{A}_T + \tilde{c}_T \xrightarrow{\mathcal{D}} Z$ as $T \to \infty$. By Proposition 2.2 $Z \overset{\mathcal{D}}{=} \mathcal{A} + s$, thus (4.27) holds. As discussed above, from (4.27) and the fact that $\tilde{c}_T \to s$ as $T \to \infty$, it follows that $\mathcal{A}_T \xrightarrow{\mathcal{D}} \mathcal{A}$.   $\square$

## 5. Extensions.

While the original motivation for our study came from the circular constraint, the proof presented extends to more general shape functions. We refrain from stating precise theorems. Still it should be instructive to the reader to see how the Brownian bridge responds to a general constraint.

Let us then substitute the circle $c_T$ by $g_T(t) = Tg(t/T)$, where $g : [-1,1] \to \mathbb{R}$, $g(-1) = g(1) = 0$, $g$ continuous, and $g \in C^2([-1,1])$ piecewise. As before we fix the reference point $\tau T$, $\tau \in (-1,1)$, and study the fluctuations away from $g_T$ for times close to $\tau T$. To first approximation the fluctuation behavior is determined by the sign of $g''(\tau)$. We list three "standard" cases, $g_c$ denoting the convex envelope of $g$.

(i) $g''(\tau) < 0$: assume that, for a $\delta > 0$, $g \in C^2$ and $g = g_c$ on $[\tau - \delta, \tau + \delta]$. If $g''(\tau) < 0$, the fluctuations are as specified in Theorem 1.1, where now $v_s = (-2g''(\tau))^{1/3}$.

(ii) $g''(\tau) = 0$: let $g$ be linear in $[t_1, t_2]$ and, for a $\delta > 0$, let $g = g_c$, $g'' < 0$, and $g \in C^2$ on $[t_1 - \delta, t_1) \cup (t_2, t_2 + \delta]$. Then the fluctuations at $t_i T$ are of order $T^\mu$, $\mu < 1/2$, and inside the interval $[t_1 T, t_2 T]$ of order $T^{1/2}$. Thus the limit process will be Brownian excursion over the interval $[t_1, t_2]$.

(iii) $g''(\tau) > 0$: let $[t_1, t_2]$ be an interval such that $t_1 < \tau < t_2$, $g(t) < g_c(t)$ for $t \in (t_1, t_2)$ and $g(t_i) = g_c(t_i)$, $i = 1, 2$. Moreover assume that for some $\delta > 0$, $g = g_c$ and is $C^2$ on $[t_1 - \delta, t_1] \cup [t_2, t_2 + \delta]$. Then in $(t_1, t_2)$ the constraint has no effect on the Brownian motion and the limit process will be a Brownian bridge over the interval $[t_1, t_2]$.

Clearly there are intermediate cases to be discussed. However, a really novel phenomenon appears if in case (i) we lift the assumption that $g$ is continuously differentiable at $\tau$. We denote the right- (left-) hand limits by $f(x^+) = \lim_{t \downarrow x} f(x)$ and $f(x^-) = \lim_{t \uparrow x} f(x)$.

(i.a) *Ridge.* Assume (i) except at $\tau$. Instead let $g''(\tau^+) < 0$, $g''(\tau^-) < 0$, and $\nu := g'(\tau^-) - g'(\tau^+) > 0$. Then the fluctuations above $g_T(\tau T)$ are of order 1 and the probability density of $X_T(\tau T)$ equals $\frac{1}{2}\nu^3 x^2 e^{-\nu x}$ as $T \to \infty$. As a consequence, (ii) and (iii) hold also if there is a ridge at $t_1$ and/or $t_2$.

(i.b) *Curvature discontinuity.* Assume (i) except at $\tau$. Instead let $g'(\tau^-) = g'(\tau^+)$ but $g''(\tau^+) \neq g''(\tau^-)$ and $g''(\tau^\pm) < 0$. In this case the fluctuations are of order $T^{1/3}$ and the limiting probability density of $X_T(\tau T)$ is, up to normalization, $\Omega(v_s(\tau^-)xT^{-1/3})\Omega(v_s(\tau^+)xT^{-1/3})$, with $v_s(\tau^\pm) = (-2g''(\tau^\pm))^{1/3}$ and $\Omega(x)$ given in (2.8).

## APPENDIX

**A.1. Properties of the Airy function and its zeros.** For the convenience of the reader we list a few properties of the Airy function needed in the main text. We follow the conventions in [1].

1. For large $z$,

(A.1) $$\mathrm{Ai}(z) \simeq \frac{1}{2\sqrt{\pi}\sqrt[4]{z}}e^{-2z^{3/2}/3}.$$

2. $\mathrm{Ai}(z) \leq 0.54$ for all $z$ and the maximum is reached at $z = -\mu \simeq -1.02$.
3. For large $k$, $\omega_k \simeq (\frac{3\pi}{2}k)^{2/3}$ and for all $k \geq 2$

(A.2) $$\omega_k - \omega_1 \geq k^{2/3}.$$

4. $|\mathrm{Ai}'(-\omega_k)| \geq \mathrm{Ai}'(-\omega_1)$ where $\omega_1 \simeq 2.34$, $\mathrm{Ai}'(-\omega_1) \simeq 0.70$.
5. For $x \in [0, -\omega_1 + \mu]$,

(A.3) $$\mathrm{Ai}(-\omega_1 + x) \geq \frac{\mathrm{Ai}(-\omega_1 + \mu)}{(-\omega_1 + \mu)}x.$$

6. For all $x \in \mathbb{R}_+$,

(A.4) $$\Omega(x) = \mathrm{Ai}(-\omega_1 + x)/\mathrm{Ai}'(-\omega_1) \leq 6e^{-x}, \qquad \Omega(x) \leq x.$$

**A.2. Leading term of the transition density.**

LEMMA A.1. *Let $\Gamma = \frac{1}{2}(t_2 - t_1)h_s$ and $y_i = v_s x_i$, $i = 1, 2$; then*

(A.5) $$\widehat{W}(x_2, t_2 | x_1, t_1) = v_s e^{-\omega_1 \Gamma}\left[\frac{\mathrm{Ai}(-\omega_1 + y_1)\mathrm{Ai}(-\omega_1 + y_2)}{\mathrm{Ai}'(-\omega_1)^2} + R_\Gamma(y_1, y_2)\right]$$

*with*

(A.6) $$|R_\Gamma(y_1, y_2)| \leq R_\Gamma^0(y_1)R_\Gamma^0(y_2),$$
$$R_\Gamma^0(y) = \min\{y\exp(-a\Gamma), \exp(-ay\Gamma^{1/3})\}$$

*for a constant $a > 0$ and $\Gamma$ large enough. Moreover, for any fixed $\Gamma > 0$, $\widehat{W}(x_2, t_2 | x_1, t_1)$ is uniformly bounded in $x_1, x_2$.*

PROOF.    Let

$$(A.7) \qquad \Phi_k(y) = e^{-(\omega_k - \omega_1)\Gamma/2} \frac{\mathrm{Ai}(-\omega_k + y)}{\mathrm{Ai}'(-\omega_k)}.$$

Then $R_\Gamma(y_1, y_2)$ is given by

$$(A.8) \qquad R_\Gamma(y_1, y_2) = \sum_{k \geq 2} \Phi_k(y_1)\Phi_k(y_2)$$

and

$$(A.9) \qquad |R_\Gamma(y_1, y_2)| \leq \sum_{k \geq 2} |\Phi_k(y_1)| \sum_{l \geq 2} |\Phi_l(y_2)|.$$

For large $k$, $\omega_k \simeq (\frac{3}{2}\pi k)^{2/3}$, and for small $k$ the exact values of $\omega_k$ are known [1], from which we deduce that $\omega_k - \omega_1 \geq \frac{1}{2}k^{2/3}$, for all $k \geq 2$. Moreover we have $|1/\mathrm{Ai}'(\omega_k)| \leq 1$ and $|\mathrm{Ai}(-\omega_k + y)| \leq |y||\mathrm{Ai}'(-\omega_k)|$. Therefore it follows that

$$(A.10) \qquad \sum_{k \geq 2} |\Phi_k(y)| \leq y \sum_{k \geq 2} e^{-k^{2/3}\Gamma/2} \leq y e^{-\Gamma/2} c_1(\Gamma)$$

with $c_1(\Gamma) = 3(\sqrt{\Gamma} + \sqrt{\pi/2})\Gamma^{-3/2}$.

This estimate is good except for very large $y$. For large $y$, the Airy function becomes of order 1 for $\omega_k \simeq y$, that is, for $k \simeq \frac{2}{3\pi}y^{3/2}$. Let $k_0(y) = y^{3/2}/10$. Then we distinguish between the cases for $k \leq k_0$ and $k \geq k_0$.

(a) $2 \leq k \leq k_0(y)$. In this case $\mathrm{Ai}(-\omega_k + y) \simeq \exp(-\frac{2}{3}(-\omega_k + y)^{3/2}) \leq \exp(-\frac{1}{3}y^{3/2})$ and, with the same estimate for the exponential term, we obtain

$$(A.11) \qquad |\Phi_k(y)| \leq \exp(-\tfrac{1}{2}k^{2/3}\Gamma)\exp(-\tfrac{1}{3}y^{3/2}).$$

(b) $k \geq k_0(y)$. For this case we use $\omega_k - \omega_1 \geq \frac{1}{2}k^{2/3}$ and $|\mathrm{Ai}(-\omega_k + y)| \leq 1$ and obtain

$$(A.12) \qquad |\Phi_k(y)| \leq \exp(-\tfrac{1}{2}k^{2/3}\Gamma).$$

Therefore for large $y$ we have

$$(A.13) \begin{aligned} \sum_{k \geq 2} |\Phi_k(y)| &= \sum_{k=2}^{k_0(y)} |\Phi_k(y)| + \sum_{k > k_0(y)} |\Phi_k(y)| \\ &\leq \sum_{k \geq 2} e^{-k^{2/3}\Gamma/2} e^{-y^{3/2}/3} + \sum_{k \geq k_0(y)} \exp(-\tfrac{1}{2}k^{2/3}\Gamma). \end{aligned}$$

The first term on the right-hand side of (A.13) is bounded by $c_1(\Gamma)\exp(-\Gamma/2 - y^{3/2}/3)$, and the second one is bounded by

$$(A.14) \qquad \int_{k_0(y)}^{\infty} dk\, e^{-k^{2/3}\Gamma/2} \leq c_2(\Gamma)e^{-\Gamma y/2}$$

with $c_2(\Gamma) = 3(\sqrt{\pi/2} + \sqrt{\Gamma y}/2)\Gamma^{-3/2}$.

If we take $\Gamma$ large, we can apply the approximation for large $y$ to the $y \geq \Gamma^{2/3}$ and use $e^{-y^{3/2}/3} \leq e^{-y\Gamma^{1/3}/3}$ to see that (A.6) holds. On the other hand, from (A.10), (A.14) and the boundedness of the ground state, it follows that $\widehat{W}(x_2, t_2 | x_1, t_1)$ is uniformly bounded in $x_1, x_2$ for any fixed $\Gamma > 0$.   □

### A.3. Estimate of the integral with error terms.

LEMMA A.2.   *Let us define*

$$I(0, \infty) = \int_0^\infty dx\, \Omega(v_j x)\Omega(v_{j+1} x)e^{-\nu x},$$

(A.15)   $$I_E(0, \infty) = \int_0^\infty dx\, \Omega(v_j x)R^0_{\Gamma_{j+1}}(v_{j+1} x)e^{-\nu x},$$

$$I_{EE}(0, \infty) = \int_0^\infty dx\, R^0_{\Gamma_j}(v_j x)R^0_{\Gamma_{j+1}}(v_{j+1} x)e^{-\nu x}.$$

*Then, if $\nu \geq 0$,*

(A.16)   $$I_E(0, \infty) \leq I(0, \infty)Ce^{-a\Gamma_{j+1}^{1/3}},$$

$$I_{EE}(0, \infty) \leq I(0, \infty)Ce^{-a(\Gamma_j^{1/3} + \Gamma_{j+1}^{1/3})}$$

*for some constant $C > 0$, assuming $\Gamma_j, \Gamma_{j+1}$ large enough.*

PROOF.   First we change variables as $y = v_j x$. Setting $\lambda = v_{j+1}/v_j$ and $\tilde{\nu} = \nu/v_j$, then

$$\tilde{I}(0, \infty) = I(0, \infty)v_j = \int_0^\infty dy\, \Omega(y)\Omega(\lambda y)e^{-\tilde{\nu}y},$$

(A.17)   $$\tilde{I}_E(0, \infty) = I_E(0, \infty)v_j = \int_0^\infty dy\, \Omega(y)R^0_{\Gamma_{j+1}}(\lambda y)e^{-\tilde{\nu}y},$$

$$\tilde{I}_{EE}(0, \infty) = I_{EE}(0, \infty)v_j = \int_0^\infty dy\, R^0_{\Gamma_j}(y)R^0_{\Gamma_{j+1}}(\lambda y)e^{-\tilde{\nu}y}.$$

To prove the lemma we have to find lower bounds for $\tilde{I}(0, \infty)$ and upper bounds for $\tilde{I}_E(0, \infty)$ and $\tilde{I}_{EE}(0, \infty)$. We use essentially (A.6), (A.3) and (A.4). First let us bound $\tilde{I}(0, \infty)$.

(a)   $\lambda \leq 1$. Let $\theta = \text{Ai}(-\omega_1 + \mu)/[(-\omega_1 + \mu)\text{Ai}'(-\omega_1)]$. Then

(A.18)   $$\tilde{I}(0, \infty) \geq \int_0^1 dy\, \theta^2 y^2 \lambda e^{-\tilde{\nu}y} = \lambda \theta^2 \kappa(\tilde{\nu}),$$

where $\kappa(\tilde{\nu}) = \int_0^1 dx\, x^2 e^{-\tilde{\nu}x}$. It is easy to see that $e^{-x} \leq 3\kappa(x)$.

(b) $\lambda \geq 1$. By the change of variable $x = \lambda y$, and then using the previous bound we obtain

$$(A.19) \qquad \tilde{I}(0,\infty) = \frac{1}{\lambda}\int_0^\infty dx\, \Omega(x)\Omega(x/\lambda)e^{-\tilde{\nu}x/\lambda} \geq \frac{1}{\lambda^2}\theta^2\kappa(\tilde{\nu}/\lambda).$$

Next we compute some upper bounds of $\tilde{I}_E(0,\infty)$.

(a) $\lambda \leq 1$.

$$(A.20) \qquad \tilde{I}_E(0,1) \leq \int_0^1 dy\, \lambda y^2 e^{-a\Gamma_{j+1}}e^{-\tilde{\nu}y} = \lambda e^{-a\Gamma_{j+1}}\kappa(\tilde{\nu})$$

and

$$(A.21) \qquad \tilde{I}_E(1,\infty) \leq \int_1^\infty dy\, \lambda y e^{-a\Gamma_{j+1}}e^{-\tilde{\nu}y}6e^{-y} \leq 6e^{-\tilde{\nu}}\lambda e^{-a\Gamma_{j+1}}.$$

(b) $\lambda \geq 1$.

$$(A.22) \qquad \tilde{I}_E(0,1/\lambda) \leq \int_0^{1/\lambda} dy\, \lambda y^2 e^{-a\Gamma_{j+1}}e^{-\tilde{\nu}y} = \frac{e^{-a\Gamma_{j+1}}}{\lambda^2}\kappa(\tilde{\nu}/\lambda)$$

and

$$(A.23) \qquad \tilde{I}_E(1/\lambda,\infty) \leq \int_{1/\lambda}^\infty dy\, y e^{-a\lambda y\Gamma_{j+1}^{1/3}}e^{-\tilde{\nu}y} \leq \frac{4}{\lambda^2}e^{-a\Gamma_{j+1}^{1/3}}e^{-\tilde{\nu}/\lambda}.$$

Putting all together, we obtain

$$(A.24) \qquad \frac{I_E(0,\infty)}{I(0,\infty)} = \frac{\tilde{I}_E(0,\infty)}{\tilde{I}(0,\infty)} \leq Ce^{-a\Gamma_{j+1}^{1/3}}$$

for all $\lambda$ with $C = 19/\theta^2$ (and $\Gamma_{j+1} \geq 1$).

Finally we bound $\tilde{I}_{EE}(0,\infty)$.

(a) $\lambda \leq 1$.

$$(A.25) \qquad \tilde{I}_{EE}(0,1) \leq \int_0^1 dy\, \lambda y^2 e^{-a(\Gamma_j+\Gamma_{j+1})}e^{-\tilde{\nu}y} = \lambda e^{-a(\Gamma_j+\Gamma_{j+1})}\kappa(\tilde{\nu})$$

and

$$(A.26) \quad \tilde{I}_{EE}(1,\infty) \leq \int_1^\infty dy\, \lambda y e^{-a\Gamma_{j+1}}e^{-ay\Gamma_j^{1/3}}e^{-\tilde{\nu}y} \leq 4e^{-\tilde{\nu}}\lambda e^{-a\Gamma_{j+1}}e^{-a\Gamma_j^{1/3}}.$$

(b) $\lambda \geq 1$. By the change of variable $x = \lambda y$ we obtain immediately $\tilde{I}_{EE}(0,1/\lambda) = \frac{1}{\lambda^2}e^{-a(\Gamma_j+\Gamma_{j+1})}\kappa(\tilde{\nu}/\lambda)$ and $\tilde{I}_{EE}(1/\lambda,\infty) \leq \frac{4}{\lambda^2}e^{-\tilde{\nu}/\lambda}e^{-a\Gamma_j}e^{-a\Gamma_{j+1}^{1/3}}$.

Putting all together, we see that for all $\lambda$

$$(A.27) \qquad \frac{I_{EE}(0,\infty)}{I(0,\infty)} = \frac{\tilde{I}_{EE}(0,\infty)}{\tilde{I}(0,\infty)} \leq Ce^{-a(\Gamma_j^{1/3}+\Gamma_{j+1}^{1/3})}. \qquad\qquad \square$$

**A.4. Monotonicity on conditioning shapes.** Let us consider a simple random walk on $\mathbb{Z}$ conditioned to come back to the origin after $2N$ steps, denoted by $\xi_N = (\xi_N(i))_{i=0}^{2N}$. Let $\Delta t = \frac{1}{2N}$, $\Delta x = \sqrt{\Delta t}$, and define $B_N(t)$ by setting $B_N(k\Delta t) = \Delta x \xi_N(k)$ for $k = 0, \ldots, 2N$, and by linear interpolation for the other values of $t \in [0,1]$. The set of possible paths $B_N$ is called $\Gamma_N$. We denote by $\mu_N$ the uniform measure on the continuous paths $B_N$.

In the sequel we consider two conditioning shapes $s_1, s_2$ such that $s_1(t) \leq s_2(t)$ for $t \in [0,1]$, $s_2(0) \leq 0$, $s_2(1) \leq 0$ and $s_2(t) < \infty$, and we denote by $\mu_N^{s_i}$ the path measure conditioned to remain above $s_i$, that is, $\mu_N^{s_i}(\cdot) = \mu_N(\cdot | B_N(t) \geq s_i(t),\ t \in [0,1])$. Let $S = C([0,1])$ be the set of bounded continuous functions from $[0,1]$ to $\mathbb{R}$ with sup norm, and define the set of increasing function by

$$(A.28) \quad \mathcal{M} = \{f \in C(S) | f(b_1) \leq f(b_2) \text{ whenever } b_1(t) \leq b_2(t) \ \forall t \in [0,1]\}.$$

PROPOSITION A.3. *If* $s_1 \leq s_2$, *then for all* $f \in \mathcal{M}$,

$$(A.29) \quad \sum_{b \in \Gamma_N} \mu_N^{s_1}(b) f(b) \leq \sum_{b \in \Gamma_N} \mu_N^{s_2}(b) f(b).$$

PROOF.   Equation (A.29) is equivalent to

$$(A.30) \begin{aligned} 0 &\leq \sum_{(b_1, b_2) \in \Gamma_N^2} \mu_N^{s_2}(b_2) f(b_2) \mu_N^{s_1}(b_1) - \sum_{(b_1, b_2) \in \Gamma_N^2} \mu_N^{s_2}(b_2) \mu_N^{s_1}(b_1) f(b_1) \\ &= \tfrac{1}{2} \sum_{(b_1, b_2) \in \Gamma_N^2} (f(b_2) - f(b_1))(\mu_N^{s_2}(b_2)\mu_N^{s_1}(b_1) - \mu_N^{s_1}(b_2)\mu_N^{s_2}(b_1)). \end{aligned}$$

Denote $\nu_N(b_1, b_2) = \mu_N^{s_2}(b_2)\mu_N^{s_1}(b_1) - \mu_N^{s_1}(b_2)\mu_N^{s_2}(b_1)$. In what follows the notation $b_1 \not\geq s_1$ means that there exists a $t$ such that $b_1(t) < s_1(t)$. Similarly, $b_1 \geq s_1$ means that $b_1(t) \geq s_1(t)$ for all $t$. For the couple $(b_1, b_2)$ there are different possibilities:

(a)  $b_1 \not\geq s_1$ and $b_2 \not\geq s_1$, then $\nu_N(b_1, b_2) = 0$.

(b)  $b_1 \geq s_1$ and $b_2 \geq s_2$, then $\nu_N(b_1, b_2) = 0$.

(c)  $b_1 \geq s_1$, $b_2 \geq s_2$, but $b_1 \not\geq s_2$, then:

    (c1)  if $b_2 \geq b_1$, then $f(b_2) - f(b_1) \geq 0$ and $\nu_N(b_1, b_2) \geq 0$ since $\mu_N^{s_2}(b_1) = 0$,

    (c2)  otherwise, $b_1$ and $b_2$ intersect above $s_2$. In this case, let $(b_1', b_2')$ be the couple of random walks defined as follows. Take a $t = t_0$ such that $b_1(t_0) < s_2(t_0)$ and set $b_1'(t_0) = b_1(t_0)$ and $b_2'(t_0) = b_2(t_0)$. Then for all other $t$ from $t_0$ to 1, $b_1'$ and $b_2'$ are defined by exchanging the paths of $b_1$ and $b_2$ when they merge and/or divide. Similarly for $t$ from $t_0$ back to 0. By the Markov property we have $\nu_N(b_1, b_2) = \nu_N(b_1', b_2')$, and the new paths satisfy

$b'_2 \geq b_1$ and $b_2 \geq b'_1$, and moreover if we apply twice the transformation we obtain the original paths. Thus, $f(b_2) + f(b'_2) - f(b_1) - f(b'_1) \geq 0$, so that the contributions in (A.30) coming from $(b_1, b_2)$ and from $(b'_1, b'_2)$ are positive.

(d) $b_2 \geq s_1$, $b_1 \geq s_2$. By symmetry the same conclusion is obtained in case (c) holds. $\square$

PROPOSITION A.4 (Invariance principle). *Let $W^0$ be the path measure of the Brownian bridge from $(0, 0)$ to $(1, 0)$. Then, as $N \to \infty$, $\mu_N \Rightarrow W^0$, that is,*

$$\text{(A.31)} \qquad \lim_{N \to \infty} \sum_{b \in \Gamma_N} \mu_N(b) f(b) = \int_S dW^0(b) f(b)$$

*for all $f \in C(S)$.*

PROPOSITION A.5. *Let $\mu^{s_i}(b) = W^0(b | b \geq s_i)$ be the path measure for the Brownian bridge conditioned to stay above $s_i$, $i = 1, 2$. We assume that $s_i$ are continuous, piecewise $C^1$, and $s_1 \leq s_2$. Then, for all $f \in \mathcal{M}$,*

$$\text{(A.32)} \qquad \int_S d\mu^{s_1}(b) f(b) \leq \int_S d\mu^{s_2}(b) f(b).$$

PROOF. Define $K(s_i)(b) = \min_{t \in [0,1]} \Theta(b(t) - s_i(t))$ with $\Theta$ the Heaviside function, and let $D_{K(s_i)}$ be the set of discontinuities of $K(s_i)$. We want to show that $\mathbb{P}_{W^0}(D_{K(s_i)}) = 0$. A path $b \notin D_{K(s_i)}$ if $\forall \varepsilon > 0$, $\exists \delta > 0$ such that $|K(s_i)(b) - K(s_i)(b')| \leq \varepsilon$, for all $b'$ satisfying $\|b' - b\|_\infty \leq \delta$. Observe that $K(s_i)(b) \in \{0, 1\}$, thus a path $b \notin D_{K(s_i)}$ if $\min_{t \in [0,1]} (b(t) - s_i(t)) \neq 0$. Therefore $b \in D_{K(s_i)}$ if $b$ touches $s_i$ but does not cross it. Now, consider a path $b$ with touches $s_i$ and let $\tau(b)$ be the first time that happens. The shape $s_i$ is continuous and piecewise $C^1$, therefore a.s. the path $b$ will cross $s_i$, thus $\mathbb{P}_{W^0}(D_{K(s_i)}) = 0$. From this follows

$$\text{(A.33)} \qquad \lim_{N \to \infty} \sum_{b \in \Gamma_N} \mu_N(b) f(b) K(s_i)(b) = \int_S dW^0(b) f(b) K(s_i)(b),$$

for all $f \in C(S)$. Since

$$\text{(A.34)} \qquad \sum_{b \in \Gamma_N} \mu_N^{s_i}(b) f(b) = \frac{\sum_{b \in \Gamma_N} \mu_N(b) f(b) K(s_i)(b)}{\sum_{b \in \Gamma_N} \mu_N(b) K(s_i)(b)},$$

(A.33) implies

$$\text{(A.35)} \qquad \lim_{N \to \infty} \sum_{b \in \Gamma_N} \mu_N^{s_i}(b) f(b) = \int_S d\mu^{s_i}(b) f(b).$$

Finally, using Proposition A.3 we conclude that (A.32) holds. $\square$

PROPOSITION A.6. *Let $\mu^{(z)}$ be the path measure for the Brownian bridge from $(0, z)$ to $(1, 0)$ conditioned to stay above a continuous piecewise $C^1$ shape $s$. If $z \geq 0$, then*

$$(A.36) \qquad \int_S d\mu^{(0)}(b) f(b) \leq \int_S d\mu^{(z)}(b) f(b)$$

*for all increasing functions $f \in \mathcal{M}$.*

PROOF. We have to show that

$$(A.37) \qquad \int_{S^2} d\mu^{(z)}(b_2) \, d\mu^{(0)}(b_1)(f(b_2) - f(b_1)) \geq 0.$$

For each couple $(b_1, b_2)$ of Brownian bridges, let $\tau(b_1, b_2) = \min_{t \in [0,1]}(b_1(t) = b_2(t))$. Define $\varphi : (b_1, b_2) \to (b'_1, b'_2)$ where $b'_i(t) = b_i(t)$ for $t \in [0, \tau(b_1, b_2)]$ and $b'_i(t) = b_{3-i}(t)$ for $t \in [\tau(b_1, b_2), 1]$, $i = 1, 2$. Obviously $\varphi(\varphi(b_1, b_2)) = (b_1, b_2)$ and by the Markov property $d\mu^{(z)}(b_2) \, d\mu^{(0)}(b_1) = d\mu^{(z)}(b'_2) \, d\mu^{(0)}(b'_1)$. By construction $b'_2 \geq b_1$, $b_2 \geq b'_1$, which implies

$$(A.38) \qquad \begin{aligned} &\int_{S^2} d\mu^{(z)}(b_2) \, d\mu^{(0)}(b_1)(f(b_2) - f(b_1)) \\ &= \tfrac{1}{2} \int_{S^2} d\mu^{(z)}(b_2) \, d\mu^{(0)}(b_1)(f(b_2) - f(b'_1) + f(b'_2) - f(b_1)) \geq 0. \end{aligned}$$

$\square$

COROLLARY A.7. *By linearity Proposition A.5 holds also if the initial and final points have a given joint density independent of the path measure.*

COROLLARY A.8. *Let $g_1, g_2$ be two probability densities such that*

$$(A.39) \qquad \int_{x \leq x_1} g_1(x) \, dx \leq \int_{x \leq x_1} g_2(x) \, dx.$$

*Denote by $\mu_x$ the path measure of Brownian motion $b(t)$ starting from $x$. Then*

$$(A.40) \qquad \begin{aligned} &\int dy \, h(y) \int dx \, g_1(x) \mu_x(f | b(1) = y) \\ &\leq \int dy \, h(y) \int dx \, g_2(x) \mu_x(f | b(1) = y) \end{aligned}$$

*for any increasing function $f \in \mathcal{M}$, where $h$ denotes the probability density of $b(1)$.*

PROOF. By linearity we need to verify the assertion only for a fixed end-point. Let $F_i(x) = \int_{y \leq x} g_i(y)\,dy$, and let $\psi_i(y) = F_i^{-1}(y)$ if $g_i(y) > 0$ and $\psi(y) = 0$ if $g_i(y) = 0$. $\psi_1(x) \leq \psi_2(x)$ for all $x$. Therefore

$$\int dx\, g_2(x)\mu_x(f|b(1) = y) = \int_0^1 dz\, \mu_{\psi_2(z)}(f|b(1) = y)$$

$$\text{(A.41)} \qquad\qquad \leq \int_0^1 dz\, \mu_{\psi_1(z)}(f|b(1) = y)$$

$$= \int dx\, g_1(x)\mu_x(f|b(1) = x). \qquad \square$$

**Acknowledgments.** We are grateful to Michael Prähofer for useful discussions and to the unknown referee for a critical reading of the previous version and insisting on convergence of path measures.

ZENTRUM MATHEMATIK
TECHNISCHE UNIVERSITÄT MÜNCHEN
D-85747 GARCHING
GERMANY
E-MAIL: ferrari@ma.tum.de
E-MAIL: spohn@ma.tum.de